\documentclass[10pt,a4paper,reqno]{amsart}

\usepackage{amsmath,amssymb,amsthm}  

\usepackage[utf8]{inputenc}

\usepackage[english]{babel}
\usepackage{url}

\usepackage[left=3.1cm, right=3.1cm]{geometry}

\def\MR#1{\href{http://www.ams.org/mathscinet-getitem?mr=#1}{MR#1}}

\usepackage{xcolor}

\definecolor{LinkColor}{rgb}{0,0,0} 
\definecolor{URLColor}{rgb}{0,.5,1}
\usepackage[colorlinks=true,linkcolor=LinkColor,citecolor=LinkColor,urlcolor=URLColor, naturalnames, hyperindex, pdfstartview=FitH, bookmarksnumbered, plainpages]{hyperref}

\newtheorem{question}{Question}

\begin{document}

\title{$3$ questions on \textsf{cut} groups}
\author{Andreas B\"achle}
\email{\href{mailto:andreas.bachle@vub.be}{andreas.bachle@vub.be}}

\maketitle

All groups considered here are finite. Rationality questions are classical and at the very heart of finite group theory. An element $x \in G$ is called \emph{rational in $G$} if it is conjugate in $G$ to all generators of $\langle x \rangle$, the cyclic subgroup generated by $x$. A group $G$ is called \emph{rational} if every element is rational in $G$. This is equivalent to the character table being a rational matrix (whence the name). The notion of rational groups was generalized by D.~Chillag and S.~Dolfi in \cite{CD}: an element $x \in G$ is called \emph{inverse semi-rational in $G$} if every generator of $\langle x \rangle$ is conjugate to $x$ or to $x^{-1}$ in $G$; a group $G$ is called \emph{inverse semi-rational} if every element is inverse semi-rational in $G$. It turned out that this has nice interpretations of different flavor. The following are equivalent for a group $G$ (see e.g.\ \cite{BCJM}):
\begin{enumerate}
 \item $G$ is inverse semi-rational.
 \item\label{cut} $\textup{Z}(\textup{U}(\mathbb{Z} G))$, the center of the unit group of the integral group ring of $G$, is finite.
 \item $K_1(\mathbb{Z} G)$, the Whitehead group of the integral group ring of $G$, is finite.
 \item For all rows of the character table of $G$, the field extension of $\mathbb{Q}$ generated by the entries of this row equals $\mathbb{Q}$ or an imaginary quadratic extension.
 \item For all columns of the character table of $G$, the field extension of $\mathbb{Q}$ generated by the entries of this column equals $\mathbb{Q}$ or an imaginary quadratic extension.
\end{enumerate}

This is a surprisingly frequently happening phenomenon: for instance, about 86.62\% of all groups up to order $512$ are inverse semi-rational (mainly due to the fact that many $2$-groups are inverse semi-rational) whereas 0.57\% of the groups of order at most 512 are rational, see \cite[Section~7]{BCJM}.  Due to the characterization \eqref{cut} above, inverse semi-rational groups are also called \emph{\textsf{cut} groups} (for \textsf{c}entral \textsf{u}nits \textsf{t}rivial, a name coined by Bakshi-Maheshwary-Passi \cite{BMP}) and for brevity this is also the term we will also use in what follows

Denote by $\mathbb{Q}(G)$ the field extension of the rationals generated by all entries of the character table of $G$. Clearly, $|\mathbb{Q}(G): \mathbb{Q}| = 1$ if and only if $G$ is rational. Is there a natural class comprising the rational groups such that the degrees of the fields $\mathbb{Q}(G)$ is uniformly bounded?

\begin{question}\label{que1} Is there $c > 0$ such that $|\mathbb{Q}(G) : \mathbb{Q}| \leqslant c$ for all \textsf{cut} groups? 
\end{question}

J.~Tent proved in \cite[Theorem~B]{Tent} that $|\mathbb{Q}(G) : \mathbb{Q}| \leqslant 2^7$ (or actually $\leqslant 2^5$) for solvable \textsf{cut} groups. There is also an affirmative answer to Question~\ref{que1} for all quasi or almost simpel groups, cf.\ \cite[Theorem~5.1]{BCJM} and S.~Trefethen's article \cite{Tre}.

On the other hand, the answer to Question~\ref{que1} is no, if one considers the slightly larger classes of semi-rational or quadratic rational groups (i.e.\ groups where one allows arbitrary quadratic extensions for each row or column of the character table, respectively) instead of \textsf{cut} groups, as can be seen from the alternating groups.

Since non-trivial rational groups have even order, the Sylow $2$-subgroups play a fundamental role in these groups. In particular it was conjectured that they should again be rational! This was refuted by I.M.~Isaacs and G.~Navarro in the article \cite{IN} providing counterexamples of order $2^9\cdot 3$, where they also proved that the Sylow $2$-subgroup of a rational group is rational again in certain classes of groups. Since every non-trivial \textsf{cut} group has an order divisbile by $2$ or $3$ \cite[Theorem~1]{BMP}, one might wonder what can be said about the corresponding Sylow subgroups. It is not hard to find examples of \textsf{cut} groups having Sylow $2$-subgroups that fail to be \textsf{cut}. However for Sylow $3$-subgroups the situation seems to be different.

\begin{question}[{\cite[Question~6.8]{BCJM}}]\label{syl} Let $G$ be a \textsf{cut} group, $P \in \operatorname{Syl}_3(G)$. Is $P$ \textsf{cut}?
\end{question}

Why might there be more hope that this question has a positive answer compared to the question on rationality of the Sylow $2$-subgroups in rational groups? One can prove the following: A $3$-element of a group $G$ is inverse semi-rational in $G$ if and only if it is inverse semi-rational in some Sylow $3$-subgroup $P$ of $G$ containing it \cite[Lemma~6.1]{BCJM}. The basic fact behind this is that the automorphism group of a cyclic $3$-group is cyclic, which is in general not the case for a cyclic $2$-group, hence the corresponding proof does not work for rationality and $2$-elements. In \cite[Section~6]{BCJM} a positive answer to Question~\ref{syl} is provided for supersolvable groups (or, more generally, for solvable groups of $3$-length $1$), Frobenius groups, for groups of small order and in several other situations. Moreover, N.~Grittini showed that the answer to Question~\ref{syl} is yes for all groups of odd order \cite[Theorem~A]{Gri}. Departing from the data in \cite[Theorem~5.1]{BCJM} and \cite{Tre} a positive answer can be obtained for all (quasi or almost) simpel groups. Note that in case Question~\ref{syl} has a positive answer, then $P/P'$ is also \textsf{cut}, hence an elementary abelian $3$-group. This is indeed always the case by a result of Isaacs-Navarro, Grittini (for solvable groups) or, in general, by \cite[Corollary~D]{NT}. 

The only primes that divide the order of a solvable rational group are $2$, $3$ and $5$ by a classical result of R.~Gow. A striking result of P.~Heged\H us asserts that $5$ plays a very special r{\^o}le: the Sylow $5$-subgroup is normal and elementary abelian in every solvable rational group \cite{Hegedus}. The only primes that divide the order of a solvable \textsf{cut} group are $2$, $3$, $5$ and $7$ \cite[Theorem~1.2]{Bac}. Yet one can construct examples of solvable \textsf{cut} groups with arbitrary large $p$-length and Sylow $p$-subgroups with arbitrary large exponent (for $p \in \{5, 7\}$ take the iterated wreath product of the normalizer of a Sylow $p$-subgroup in $S_p$, the symmetric group of degree $p$). But is at least the Hall $\{5,7\}$-subgroup of each Fitting layer nice? As usual, $\operatorname{O}_p(G)$ denotes the largest normal $p$-subgroup of $G$.

\begin{question} Let $G$ be a solvable \textsf{cut} group. Is it true that $\exp \operatorname{O}_p(G) \mid p$ for $p \in \{5, 7\}$? \end{question}


\providecommand{\bysame}{\leavevmode\hbox to3em{\hrulefill}\thinspace}
\providecommand{\MR}{\relax\ifhmode\unskip\space\fi MR }
\providecommand{\MRhref}[2]{%
  \href{http://www.ams.org/mathscinet-getitem?mr=#1}{#2}
}
\providecommand{\href}[2]{#2}

\end{document}